\documentclass[a4paper, 10pt]{article}
\usepackage{amsmath}
\usepackage{amsfonts}
\usepackage{amsthm}
\usepackage[left=2.7cm,right=2.7cm,top=3.5cm,bottom=3cm]{geometry}
\usepackage{amsmath,amssymb,amsthm,amscd,amsfonts}
\usepackage{latexsym}
\usepackage{graphicx}
\usepackage[all]{xy}
\usepackage{xcolor}
\usepackage[OT2,T1]{fontenc}
\usepackage{hyperref}
\usepackage{mathtools}
\usepackage{tikz}
\usetikzlibrary{arrows}
\usepackage{float}
\usepackage{caption}
\usepackage[new]{old-arrows}

\theoremstyle{plain}
\newtheorem{theorem}{Theorem}[section]
\newtheorem{lemma}[theorem]{Lemma}
\newtheorem{proposition}[theorem]{Proposition}

\newtheorem{Problem}[theorem]{Problem}

\theoremstyle{definition}
\newtheorem{definition}[theorem]{Definition}

\theoremstyle{remark}
\newtheorem{remark}[theorem]{Remark}

\DeclareMathAlphabet{\mathpzc}{OT1}{pzc}{m}{it}
\DeclareSymbolFont{cyrletters}{OT2}{wncyr}{m}{n}
\DeclareMathSymbol{\Sha}{\mathalpha}{cyrletters}{"58}

\begin{document}
\thispagestyle{empty}

\newcommand{\Z}{{\mathbb Z}}
\newcommand{\Q}{\mathbb Q}
\newcommand{\FF}{{\mathbb F}}
\newcommand{\E}{{\mathcal{E}}}
\newcommand{\mcF}{{\mathcal{F}}}
\newcommand{\A}{\mathcal{A}}
\newcommand{\G}{\mathbb{G}}
\newcommand{\loc}{{\rm{loc}}}
\newcommand{\chara}{\rm{char}}
\newcommand{\ds}{\displaystyle}
\newcommand{\la}{\langle}
\newcommand{\ra}{\rangle}
\newcommand{\z}{{\zeta}}
\newcommand{\ov}{\overline}
\newcommand{\wt}{\widetilde}
\newcommand{\Or}{\mathcal{O}}
\newcommand{\X}{\mathcal{X}}

\newcommand{\Hil}{\mathcal{H}}
\newcommand{\End}{\rm{End}}
\newcommand{\Aut}{\rm{Aut}}

\newcommand{\Gal}{{\rm Gal}}
\newcommand{\Id}{{\rm Id}}
\newcommand{\GL}{{\rm GL}}
\newcommand{\SL}{{\rm SL}}
\newcommand{\Hh}{{\rm H}}
\newcommand{\Mat}{{\rm Mat}}
\newcommand{\disc}{{\rm disc}}
\newcommand{\imm}{{\rm Im}}
\newcommand{\ind}{{\rm Ind}}
\newcommand{\res}{{\rm res}}
\newcommand{\diag}{{\rm diag}}
	
\renewcommand{\leq}{\leqslant}
\renewcommand{\geq}{\geqslant}
\newcommand{\modn}{{\rm{mod} \hspace{0.1cm} }}
\newcommand{\nc}{\normalcolor}
\newcommand{\teal}{\color{teal}}
\newcommand{\rc}[1]{{\color{black}{#1}}}
\newcommand{\ja}[1]{{\color{black}{#1}}}
\newcommand{\lp}[1]{{\color{black}{#1}}}

\title{Local-global divisibility on algebraic tori}
\author{Jessica Alessandr\`i, Rocco Chiriv\`i, Laura Paladino}
\date{}
	
\maketitle

\renewcommand{\thefootnote}{\arabic{footnote}}
\setcounter{footnote}{0}
	
\begin{abstract}
    We give a complete answer to the local-global divisibility problem for algebraic tori. In particular, we prove that given an odd prime $p$, if $T$ is an algebraic torus of dimension $r< p-1$ defined over a number field $k$, then the local-global divisibility by any power $p^n$ holds for $T(k)$. We also show that this bound on the dimension is best possible, by providing a counterexample for every dimension $r \geq p-1$. Finally, we prove that under certain hypotheses on the number field generated by the coordinates of the $p^n$--torsion points of $T$, the local-global divisibility still holds for tori of dimension less than $3(p-1)$.
\end{abstract}
\medskip

\medskip

\textbf{Mathematics~Subject~Classification~(2020)}: 11E72, 14G05, 11G35.

\section{Introduction}

Let $k$ be a number field and let $\mathcal G$ be a commutative algebraic group defined over $k$. We denote by $M_k$ the set of places of $k$ and by $k_v$ the completion of $k$ at $v$.
In 2001, Dvornicich and Zannier, motivated by a particular case of the Hasse principle for binary quadratic forms, stated the following problem, which is known as the \emph{Local-global divisibility problem} (see \cite{DZ1}).
\begin{Problem}[Dvornicich and Zannier, \cite{DZ1}]\label{prob}
	Let $q$ be a fixed positive integer. If we assume that the point $P \in \mathcal G(k)$ has the following property: for all but finitely many $v\in M_k$ there exists $D_v \in \mathcal G(k_v)$ such that $P = qD_v$; can we conclude that there exists $D \in \mathcal G(k)$ such that $P = qD$?
\end{Problem}
Clearly it is sufficient to answer the problem when $q$ is a power of a prime. 
The classical case of the multiplicative group $\mathcal G = \G_m$ has a complete answer: positive for $q$ \ja{not} divisible \ja{by $8$} (see for example \cite{DZ1,AT}) and negative for $q$ divisible by \ja{$8$} (see \cite{Trost}).
For a general commutative algebraic group, Dvornicich and Zannier gave in \cite{DZ1} a cohomological interpretation of the problem and sufficient conditions to answer the question (see also \cite{DZ3}). In particular, they showed that such an answer is linked to the triviality of a subgroup of $\Hh^1(\Gal\left(k\left(\mathcal G[q]\right)/k\right), \mathcal G[q])$, where $k\left(\mathcal G[q]\right)$ is the field obtained by adjoining the coordinates of the $q$--torsion points of $\mathcal G$ to $k$, called the \emph{first local cohomology group} (see Section \ref{prelim} for further details).

By using these tools, it was possible to give criteria to answer the local-global divisibility problem for several algebraic groups.
 In the case of elliptic curves the problem was long studied. 
The answer 
is affirmative for $q=2,3$ and for all powers $q=p^n$, with $p> C([k:\Q])$, where
$C([k:\Q])=3$, if $k=\Q$, and $C([k:\Q])=\left(3^{\frac{[k:\Q]}{2}}+1\right)^2$, if $k\neq \Q$ (see the recent survey \cite{DP} for further
 details; see also \cite{Creu3}). Instead, for  $q=p^n$, with $p=2,3$ and $n\geq 2$, there are known counterexamples over $\Q$ and over $\Q(\z_3)$ (when
 $p=3$) \cite{Pal3, Creu2}.
For each number field $k$ linearly disjoint over $\Q$ (\ja{resp.\ }over $\Q(\z_3)$) from the $p^n$--division field of $\E$ over $\Q$ (\ja{resp.\ }over $\Q(\z_3)$) these counterexamples \ja{also give} counterexamples in a finite extension of $k$ (see Remark
\ref{extended_counterexamples} for further details).
While there are no explicit counterexamples for $p\geq 5$ over a number field $k$, in \cite{Ran} necessary conditions on $\Gal(k\left(\E[p]\right)/k)$ that have to be satisfied in order to have local-global divisibility by $p^n$ for a prime $p\geq 5$ are given. 
In addition, for elliptic curves an effective version of the hypotheses of Problem \ref{prob} is produced in \cite{DP2}. 
For principally polarized abelian surfaces in \cite{GRab} Gillibert and Ranieri proved sufficient conditions for the local-global divisibility by any prime power $p^n$, while in \cite{GRgl2} they \ja{generalized} these conditions in order to answer the case of $\GL_2$--type varieties (see also \cite{GRgl}).
Furthermore, in \cite{Pal4} the third author produced conditions for the local-global $p$--divisibility for a general commutative algebraic group. 
In the case of abelian varieties, the problem is also linked to a classical question posed by Cassels in 1962 on the $p$--divisibility of the Tate-Shafarevich group (see \cite{Cas62, CS, Creu, DP}). 
	
\bigskip
In this work we focus on algebraic tori. As mentioned above, for the one-dimensional split torus $\G_m$ we have a complete answer. Notice that the negative answer for $q=2^n$, with $n \geq 3$, implies that one can find counterexamples in every dimension, just by taking  direct products of copies of $\G_m$. In \cite{DZ1} Dvornicich and Zannier proved that the local-global divisibility by a prime $p$ holds for tori of dimension $r \leq \max\{3,2(p-1)\}$, but fails for a torus of dimension $r = p^4-p^2+1$. In this last case they produced a $k$--rational point which is locally $p$--divisible for \lp{all but finitely many places $v\in M_k$, but not globally $p$--divisible.}
In \cite{Ill} Illengo improved the condition $r\leq 2(p-1)$ given by Dvornicich and Zannier with the weaker one $r<3(p-1)$. He also proved that this bound is best possible, by building an example with $r=3(p-1)$ for which the local-global divisibility by $p$ fails.

In this work we prove that the local-global divisibility by any odd power $p^n$, with $n\geq 1$, holds for an algebraic torus of dimension $r<p-1$. 
In particular, we prove that if $r<p-1$, then we have an affirmative answer to Problem \ref{prob}, while for $r\geq p-1$ the local-global divisibility by $p^n$, with $n \geq 2$, is no longer assured. For the latter, we construct a counterexample of dimension $p-1$ for which the local-global divisibility by every $p^n$ with $n \geq 2$ does not hold. We remark that, starting from this construction, one can build a counterexample for every dimension $r \geq p-1$, by taking the product of the torus that we build in Lemma \ref{costr toro} with the split torus of dimension $r-(p-1)$. 
We summarize these results in the following theorem, that we prove in Section \ref{sec:dim1}:

\begin{theorem}\label{thm:1}
	Let $p$ be an odd prime.
	\begin{itemize}
		\item[(a)] Let $k$ be a number field and let $T$ be a torus defined over $k$. If $T$ has dimension less than $p-1$, then the local-global divisibility by $p^n$ holds for $T(k)$, for every $n\geq 1$.
		\item[(b)] For \ja{every $n \geq 2$ and for every} $r\geq p-1$ there exists a torus defined over $k = \Q(\z_p)$ of dimension $r$ and a finite extension $L/k$ such that the local-global divisibility by $p^n$ does not hold for $T(L)$.
	\end{itemize}
\end{theorem}
Nevertheless, under certain conditions on the base field $k$, \ja{satisfied e.g.\ by} adjoining a primitive $p^n$--th root of unity to $k$, we can say something more when $p-1 \leq \dim(T)<3(p-1)$. Let $k(T[p^n])$ be the field obtained by adjoining the coordinates of the $p^n$--torsion points of $T$ to $k$. In Section \ref{sec:dim2} we prove the following.

\begin{theorem}\label{thm:2}
	Let $p$ be an odd prime and let $n \geq 1$ be an integer. Suppose that $T$ is a torus defined over $k$ with $p-1\leq \dim(T) < 3(p-1)$ and $p$ does not divide the degree $[k\left(T[p^n]\right) \cap k\left(\z_{p^n}\right) : k]$, where $\z_{p^n}$ is a $p^n$--th root of unity. Then the local-global divisibility by $p^n$ holds for $T(k)$.
\end{theorem}

We remark that, since the proof of this theorem is done by induction on the powers of $p$, we need the condition $\dim(T) < 3(p-1)$ for the base of the induction.

\bigskip
\noindent{\bf Acknowledgements.} The authors warmly thank Davide Lombardo for his precious help during the preparation of this manuscript and Umberto Zannier for useful discussions. The first and third authors are members of INdAM-GNSAGA. Part of this research was done while the first author was visiting the third author; the first author is grateful to the Department of Mathematics and Computer Science of the University of Calabria for the kind hospitality and to the University of L'Aquila for the relative funding. We also thank the referee for the helpful comments and suggestions that allowed us to improve the results of Theorem \ref{thm:1} and lighten the proof of Lemma \ref{costr toro}.

\section{Preliminary Results}\label{prelim}
As mentioned above, a useful method introduced by Dvornicich and Zannier in \cite{DZ1} in dealing with this problem is to translate it into cohomological terms. Let us recall some definitions and results.

Let $\mathcal G$ be a commutative algebraic group, defined over a number field $k$. Given a positive integer $q$, we denote by $\mathcal G[q]$ the set of $q$--torsion points of $\mathcal G(\ov k)$. It is isomorphic to $\left(\Z/q\Z\right)^\ell$, for some $\ell$ depending only on $\mathcal G$ (see for example \cite[Section 2]{DZ1}). 

Let $K := k(\mathcal G[q])$ be the number field generated by adjoining to $k$ the coordinates of the $q$--torsion points of $\mathcal G$. It is a Galois extension of $k$ and we denote by $G$ its Galois group. Let $\Sigma$ be the set of places $v$ unramified in $K$ and let $G_v = \Gal(K_w/k_v)$, where $w$ is a place of $K$ extending $v \in \Sigma$.

In \cite{DZ1}, the authors introduce a subgroup of $\Hh^1(G,\mathcal G[q])$ called the \emph{first local cohomology group}:
\begin{equation}\label{def h1loc}
	\Hh_{\loc}^1(G,\mathcal G[q]) = \bigcap_{v \in \Sigma} \ker \left(\Hh^1(G,\mathcal G[q]) \xrightarrow{\res_v} \Hh^1(G_v, \mathcal G[q])\right).
\end{equation}
By \v{C}ebotarev Density Theorem, the Galois groups $G_v$ run over all cyclic subgroups of $G$, as $v$ varies in $\Sigma$.  Therefore a cocycle $\{Z_\sigma\}_{\sigma \in G}$ in $\Hh^1(G,\mathcal G[q])$ is an element in $\Hh_{\loc}^1(G,\mathcal G[q])$ if and only if for every $\sigma \in G$ there exists $W_{\sigma} \in \mathcal G[q]$ such that $Z_\sigma = (\sigma-1)W_\sigma$ (see \cite{DZ1} and \cite{DP} for further details). This justifies the following, more general, definition.

\begin{definition}[Dvornicich and Zannier \cite{DZ1}]
	Let $G$ be a group and let $M$ be a $G$--module. We say that a cocycle $ \{Z_g\}_{g \in G}$ of $G$ with values in $M$ \emph{satisfies the local conditions} if there exist $W_g \in M$ such that $Z_g= (g-1)W_g$ for all $g\in G$. We denote by $\Hh^1_\loc(G,M)$ the subgroup of $\Hh^1(G,M)$ of the classes of these cocycles.
\end{definition}
	
As mentioned, it is sufficient to consider the case when $q$ is a prime power, thus from now on let $q= p^n$, for some prime number $p$ and some positive integer $n$. By \cite[Proposition 2.1]{DZ1} the triviality of the first local cohomology group gives a sufficient condition for an affirmative answer to Problem \ref{prob}.
\begin{proposition}[Dvornicich and Zannier \cite{DZ1}]\label{suffcond}
	Assume that $\Hh^1_{\loc}(\Gal(K/k), \mathcal G[q]) = 0$. Let $P \in \mathcal G(k)$	be a rational point with the following property: for all but finitely many primes $v$ of $k$, there exists $D_v \in \mathcal G(k_v)$ such that $P = qD_v$. Then there exists $D \in \mathcal G(k)$ such that $P = qD$.
\end{proposition}

Furthermore, by \cite[Theorem 2.5]{DZ1} it is enough to prove that $\Hh^1_{\loc}(G_p, \mathcal G[q]) = 0$ for a $p$--Sylow subgroup $G_p$ of $G$. Moreover, the converse of Proposition \ref{suffcond} is also true over a finite extension of $k$, as the following theorem by Dvornicich and Zannier shows.

\begin{theorem}[{\cite[Theorem 3]{DZ3}}]\label{neccond}
	Suppose that $\Hh^1_\loc(G, \mathcal G[q])$ is not trivial. Then there exists a number field $L$ such that $L \cap K = k$ and a point $P \in \mathcal G(L)$ which is divisible by $q$ in $\mathcal G(L_w)$ for all places $w$ of $L$ but is not divisible by $q$ in $\mathcal G(L)$.
\end{theorem}
	
We remark that if in \eqref{def h1loc} we take all valuations, instead of almost all, then we get a group isomorphic to the Tate-Shafarevich group $\Sha(k,\mathcal G[q])$. Thus the vanishing of $\Hh^1_\loc(G,\mathcal{G}[q])$ implies the vanishing of $\Sha(k,\mathcal G[q])$, which is a sufficient condition to give an affirmative answer to the local-global divisibility problem in the case when $v$ runs over all valuations of $k$ (see for instance \cite{Creu2}).
\medskip
    
Let us now introduce some notation for algebraic tori that we will use in the next sections. We adopt the same notation as in \cite[Section 4]{DZ1}. Let $T$ be an algebraic torus defined over a number field $k$, of dimension $r$.
There exists an isomorphism of algebraic groups (defined over $\ov k$) $\phi:T \longrightarrow \G^r_m$. Let $G_k$ be the absolute Galois group of $k$.
For $\sigma \in G_k$, we denote by $\phi^{\sigma}$ the twist of $\phi$ by $\sigma$, i.e. $\phi^\sigma = \sigma\circ \phi \circ \sigma^{-1}$, and we consider the following map 
\[\begin{matrix}
	\psi: & G_k & \longrightarrow & \Aut(\G^r_m) & \simeq \GL_r(\Z)\\
		& \sigma & \longmapsto & \phi\circ (\phi^\sigma)^{-1}.\\
\end{matrix}
\]
Notice that $\psi$ is a 1--cocycle, but since the action of $G_k$ on $\G^r_m$ is trivial, $\psi$ is actually a group homomorphism.

\noindent The isomorphism $\phi$ is defined over some number field, thus $\ker\psi$ has finite index in $G_k$ and we can identify $\Delta:= \psi(G_k)$ with a finite subgroup of $\GL_r(\Z)$. We denote by $L$ the field fixed by $\ker\psi$; it is a normal extension of $k$ and $\Delta \simeq \Gal(L/k)$. The field $L$ is also the splitting field of the torus $T$.
Let $\z :=\z_{p^n}$ be a primitive $p^n$--th root of unity and let $\chi$ be the cyclotomic character
\[
\begin{matrix}
	\chi: & G_k & \longrightarrow & \left(\Z/p^n\Z\right)^{\times}\\
		 & \sigma & \longmapsto & j_\sigma,\\
\end{matrix}
\]
where $j_\sigma$ is such that $\sigma(\z) = \z^{j_\sigma}$.
Let $T[p^n]$ be the group of the $p^n$--torsion points of $T$. We have $T[p^n] = T(\ov k)[p^n] \simeq \left\{ \left(\z^{j_1},\dots,\z^{j_r}\right) \in \left(\ov k^{\times}\right)^r \mid j_h \in \Z/p^n\Z \right\}$ and we fix the following isomorphism
\begin{equation}\label{iso}
\begin{matrix}
	T[p^n] & \longrightarrow & \left(\Z/p^n\Z\right)^r\\
	(\z^{j_1},\dots,\z^{j_r}) & \longmapsto & (j_1,\dots,j_r).\\
\end{matrix}
\end{equation}
By this isomorphism, the natural action of $G_k$ on $T[p^n]$ induces the following action on $\left(\Z/p^n\Z\right)^r$: $\sigma \cdot v = j_\sigma \wt{\psi(\sigma)} v$ for all $v \in \left(\Z/p^n\Z\right)^r$, where the tilde denotes the reduction $\modn p^n$. Therefore we have the homomorphism
\begin{equation}\label{xi}
    \begin{matrix}
	   \xi: & G_k & \longrightarrow & \GL_r(\Z/p^n\Z)\\
		& \sigma & \longmapsto & j_\sigma\wt{\psi(\sigma)}.\\
	\end{matrix}
\end{equation}

It is easy to check that the field fixed by $\ker\xi$ is $K := k\left(T[p^n]\right)$, thus the image $G$ of $\xi$ in $\GL_r(\Z/p^n\Z)$ is a finite subgroup isomorphic to $G_k/\ker\xi \simeq \Gal(K/k)$. Let $G_{k\left(\z\right)}$ be the subgroup $\Gal\left(\ov{k}/k\left(\z\right)\right)$ of $G_k$. We have $L = \ov{k}^{\ker\psi}$ and $K = \ov{k}^{\ker\xi}$, since $\ker\xi \supseteq \ker\psi\cap G_{k\left(\z\right)}$ we also get $K \subseteq L\left(\z\right)$.
The kernel of the restriction of $\xi$ to $G_{k\left(\z\right)}$ is contained both in $\ker\xi$ and in $\ker \wt\psi$. It follows that the image of $G_{k\left(\z\right)}$ via $\xi$ is a normal subgroup $G'$ of $G$, which is also a normal subgroup of the reduction $\wt\Delta$ modulo $p^n$ of $\Delta$. In particular, we have the following tower of extensions:

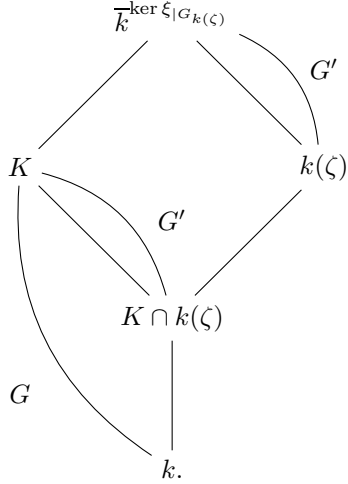
\begin{figure}[H]
	\centering
		\begin{tikzpicture}
			\node (a) at (0,0) {$\overline k ^{\ker \xi_{|G_{k(\z)}}}$};
			\node (b) at (-2,-2) {$K$};
			\node (c) at (2,-2) {$k(\z)$};
			\node (d) at (0,-4) {$K\cap k(\z)$};
			\node (e) at (0,-6) {$k$.};
			\draw (a) -- (b);
			\draw (a) -- (c);
			\draw (b) -- (d);
			\draw (c) -- (d);
			\draw (d) -- (e);
			\node (f) at (2,-0.7) {$G'$};
			\node (g) at (0,-2.7) {$G'$};
			\node (h) at (-2,-5) {$G$};
			\draw[bend left] (a) edge (c);
			\draw[bend left] (b) edge (d);
			\draw[bend right] (b) edge (e);
		\end{tikzpicture}
		\caption*{}
\end{figure}

We have obtained that $G$ and $\wt\Delta$ have a common normal subgroup $G'$, with $[G:G'] \mid p^{n-1}(p-1)$ and $[\wt\Delta:G']\mid p^{n-1}(p-1)$. In \cite{DZ1} the authors can easily conclude that $G$ and $\wt\Delta$ also have the same $p$--Sylow subgroups (since in their situation both $[G:G']$ and $[\wt\Delta:G']$ are coprime with $p$) and thus they reduce themselves to study only the $p$--Sylow subgroups in $\wt\Delta$. Instead, in our general setting for the divisibility by $p^n$ problem, we have to distinguish two cases: either $p \mid [G:G']$ or $p \nmid [G:G']$.
		

The following theorem (see \cite[Theorem 6.1.16]{Scott}) will be a precious tool in proving part (a) of Theorem \ref{thm:1}.
\begin{theorem}\label{sylow quot}
    A subgroup of a quotient is a $p$--Sylow subgroup if and only if it is the image through the canonical projection homomorphism of a $p$--Sylow subgroup.
\end{theorem}


In the proof of Theorem \ref{thm:2} we will also need the following well known result (see p.197 of \cite{Mink}) whose proof we include for the reader's convenience (for a more general result see \cite{Serre}).
\begin{lemma}\label{lemma:inj}
	Let $p$ be an odd prime and let $\pi :\GL_r(\Z)\longrightarrow \GL_r(\Z/p\Z)$ be the reduction modulo $p$. Then $\pi$ is injective on finite subgroups of $\GL_r(\Z)$.
\end{lemma}
\begin{proof}
	It is enough to show that if $A\in \ker\pi$ and $A$ has finite order $m$, then $A=\Id$. Suppose that instead $A\neq \Id$. If $p\nmid m$, then write $A = 1 + p^k B$, with $k\geq 1$ and $B\in \Mat_r(\Z)$ such that $p$ does not divide at least one of the entries of $B$. We have
	\begin{equation*}
		1 = A^m = \sum_{j=0}^m \binom{m}{j}p^{jk}B^j = 1 + mp^kB + \sum_{j=2}^m \binom{m}{j}p^{kj}B^j.
	\end{equation*}
	Thus $mp^kB = -p^{k+1}C$, for some $C \in \Mat_r(\Z)$ and so $mB = -pC$, which is a contradiction. On the other hand, if $p \mid m$, then $\ov A := A^{m/p}$ has order $p$ and lies in $\ker \pi$. We have $\ov A = 1 + p^k\ov B$ with $k\geq 1$ and $\ov B\in \Mat_r(\Z)$ such that $p$ does not divide at least one of the entries of $\ov B$. We have
	\begin{equation*}
		1 = \ov A^p = \sum_{j=0}^p \binom{p}{j}p^{jk}\ov B^j = 1 + p^{k+1}\ov B + \sum_{j=2}^{p-1}{\binom{p}{j}p^{kj}\ov B^j} + p^{pk}\ov B^p.
	\end{equation*}
	Since $p\neq 2$ and $k \geq 1$, we get $p^{k+1}\ov B = -p^{k+2} \ov C$, for some $\ov C \in \Mat_r(\Z)$, so $B = -p\ov C$ and we have again a contradiction.

	\qedhere
\end{proof}	

\noindent Finally, we state the following Lemma by Illengo; this is a key result for the local-global divisibility in algebraic tori and we will use it in the proof of Theorem \ref{thm:1}.
\begin{lemma}[{\cite[Lemma 4]{Ill}}]\label{illengo}
	Let $p$ be a prime and let $\Gamma$ be a $p$--group of matrices in $\SL_r(\Q)$. If $r<p(p-1)$ then $\Gamma$ is isomorphic to $(\Z/p\Z)^b$, for some $b \leq r/(p-1)$.
\end{lemma}
	
\section{Proof of Theorem \ref{thm:1}}\label{sec:dim1}
In this section we prove Theorem \ref{thm:1}. With the notation as above, we show that:
\begin{itemize}
	\item[(a)] let $p$ be an odd prime number and let $n\geq 1$ be an integer, then for every algebraic torus $T$ of dimension $r<p-1$ we have $\Hh^1_\loc(G_p, T[p^n]) = 0$, where $G_p$ is a $p$--Sylow subgroup of $G$. Hence by \cite[\ja{Proposition 2.5}]{DZ1} Problem \ref{prob} has affirmative answer;
	\item[(b)] for every odd prime number $p$ and every positive integer $n\geq 2$, there exists a torus $T$ defined over $\Q(\z_p)$ of dimension $p-1$ such that $\Hh^1_{\loc}\left(G, T[p^n]\right) \neq 0$. Thus by Theorem \ref{neccond}, there exists a number field in which the local-global divisibility by $p^n$ with $n \geq 2$ does not hold in $T$.
\end{itemize}
The counterexample in (b) shows that the bound on the dimension of $T$ is best possible.

Since point (b) requires more effort we start by showing a few results that we will use for its proof. The first step is building a torus such that $G$ is isomorphic to $\Z/p\Z \times \Z/p^{n-1}\Z$.
\begin{lemma}\label{costr toro}
	Let $p$ be an odd prime \ja{and let $n \geq 2$}. There exists an algebraic torus $T$ of dimension $r = p-1$ defined over $k =\Q(\z_p)$ such that $G$ is isomorphic to $\Z/p\Z\times \Z/p^{n-1}\Z$. In particular, $G \subseteq \GL_r(\Z/p^n\Z)$ is generated by 
	\[
	   \gamma_1 =\begin{pmatrix}
			0 &  & & & -1\\
			1 & 0 &  &  & -1\\
			 & \ddots & \ddots & & \vdots\\
			& & 1 & 0 & -1\\
			 &   & & 1 & -1\\
		\end{pmatrix}
		\quad \text{and} \quad \gamma_2 =
		\begin{pmatrix}
			p+1 & & \\
		     & \ddots &\\
			 && p+1\\
		\end{pmatrix}. 
	\]
\end{lemma}
\begin{proof}
	Let $L$ be a Kummer extension of $k$, such that $[L:k]=p$ and a prime other than $p$ ramifies, e.g.\ take $L = k(\sqrt[p]{2})$. Then $L \cap \Q(\z_{p^n})= \Q$ and $L/k$ is a cyclic extension of degree $p$.
	Let $\sigma$ be a generator of $\Gal(L/k)$.
	Consider the split torus $\G_m = \G_{m,L}$ defined over $L$. We denote the group $\G_m[p^n]$ by $V\simeq \Z/p^n\Z$ using additive notation. The Galois action on $V$ is given by the cyclotomic character $\chi: G_k \longrightarrow (\Z/p^n\Z)^{\times}$. 
	Let $X = R_{L/k}\G_m$ be the Weil restriction of $\G_m$. It is an algebraic torus defined over $k$ of dimension $p$, split over $L$. The group $X[p^n]$ of the $p^n$--torsion points of $X$ is a free $\Z/p^n\Z$--module of rank $p$ and by properties of the Weil restriction we have that 
    \[X[p^n] \simeq \prod_{j=0}^{p-1} V_j,\]
    where $V_j$ is an isomorphic copy of $V$, for every $j$. Let us choose a lifting $\overline{\sigma} \in G_k$ of $\sigma$. For every $\tau \in G_L = \Gal(\overline{k}/L)$ and $0 \leq h\leq p-1$, the action of an element $\gamma = \tau\overline{\sigma}^h \in G_k$ on $X[p^n]$ is given by  
    \begin{align}\label{xi action}
		\gamma \cdot (x_0,x_1\dots,x_{p-1}) & = \tau\overline{\sigma}^h \cdot (x_0,x_1\dots,x_{p-1}) \nonumber \\
		& = \chi(\gamma) (x_{-h},x_{-h+1},\dots,x_{-h-1}),
	\end{align}
    where we are considering the indices of the coordinates modulo $p$.

	Now let $T$ be the norm $1$ subtorus of $X$, that is $T = R^{(1)}_{L/k} \G_m = \ker (R_{L/k}\G_m \stackrel{N_{L/k}}{\longrightarrow} \G_m)$, the kernel of the (generalized) norm map on $X$. It is an algebraic torus over $k$ of dimension $r = p-1$, split over $L$. We are going to show that $k\left(T[p^n]\right)= L\left(\z_{p^n}\right)$.\\
	Through the isomorphism \eqref{iso} applied to $X[p^n]$, we can regard $T[p^n]$ as the submodule $W$ of $X[p^n]$ of \ja{those} vectors $(x_0,x_1\dots,x_{p-1})$ such that the sum of all coordinates is equal to zero (we are using the additive notation here).
	The Galois action on $T[p^n]$ is given by $\xi: G_k \longrightarrow \Aut(W)$, i.e.\ by the action on the points of $X[p^n]$ that lie in $W$ (see \eqref{xi}). We have that $k\left(T[p^n]\right) = \ov k^{\ker\xi}$.
	Thus, in order to determine this field, we need to find the elements of $G_k$ that act trivially on $W$.\\
	Since $p\geq 3$, by \eqref{xi action} we see that $\gamma = \tau\overline{\sigma}^h$ \ja{ fixes every $(x_0,x_1\dots,x_{p-1})$ in $W$ if and only if $h = 0$ and $\chi(\gamma)=\chi(\tau) = 1$}, i.e.\ if and only if \ja{$\gamma \in G_L \cap G_{\Q(\z_{p^n})} = G_{L(\z_{p^n})}$}.
	We therefore conclude that $k\left(T[p^n]\right) = L(\z_{p^n})$ as claimed.
	
    \ja{Note that} $L \cap \Q(\z_{p^n}) = \Q(\z_p) = k$, \ja{so} the extension $k\left(T[p^n]\right)/k$ has Galois group $\Gal\left(k\left(T[p^n]\right)/k\right)$ and
    \begin{equation}\label{eq:isogal}
    \Gal\left(k\left(T[p^n]\right)/k\right)\ni\varphi\longmapsto(\varphi_{|L}, \varphi_{|k(\zeta_{p^n})})\in\Gal\left(L/k\right) \times \Gal\left(k(\z_{p^n})/k\right)
    \end{equation}
    is an isomorphism.
    Further, it is clear that the last group is isomorphic to $\Z/p\Z \times \Z/p^{n-1}\Z$.
    
    \ja{\noindent Let $\eta\in \Gal\left(k(\z_{p^n})/k\right)$ be the automorphism sending $\z_{p^n}$ to $\z_{p^n}^{p+1}$; the two elements $\sigma$ and $\eta$ are generators of $\Gal\left(L/k\right) \times \Gal\left(k(\z_{p^n})/k\right)$.}
    As noticed in the previous section, the group $\Gal\left(k\left(T[p^n]\right)/k\right)$ is isomorphic to $G = \xi(G_k) \subseteq \GL_{p-1}(\Z/p^n\Z)$. So we want to represent $\sigma$ and $\eta$ as matrices in $\GL_{p-1}(\Z/p^n\Z)$.
    We can choose the lifting $\overline\sigma \in G_k$ of $\sigma$ such that, when restricted to $k\left(T[p^n]\right)$, it corresponds to the pair $(\sigma, 1)$ in the isomorphism \eqref{eq:isogal}; in particular $\chi(\overline\sigma) = 1$.
    With respect to the the basis $v = (1,-1, 0, \ldots, 0,0),\,\overline\sigma(v),\,\overline\sigma^2(v),\ldots,\overline\sigma^{p-2}(v)$ of $W$, we can write the matrix  $\gamma_1 = \xi(\overline\sigma)$ corresponding to $\sigma$, as
    \[
        \gamma_1 =\begin{pmatrix}
    			0 &  & & & -1\\
    			1 & 0 &  &  & -1\\
    			 & \ddots & \ddots & & \vdots\\
    			& & 1 & 0 & -1\\
    			 &   & & 1 & -1\\
    		\end{pmatrix}.
    \]
    To conclude the proof, we observe that, with a similar reasoning, we can lift $\eta$ to $\overline \eta\in G_k$ in such a way that, when restricted to $k\left(T[p^n]\right)$, the element $\overline\eta$ corresponds to the pair $(1,\eta)$ and, clearly, $\chi(\overline\eta) = p+1$. Hence the action of $\overline\eta$ on $W$ is just by multiplication by $p+1$; so if $\gamma_2 = \xi(\overline\eta)$ we have
    \[\gamma_2 =
		\begin{pmatrix}
			p+1 & & \\
		     & \ddots &\\
			 && p+1\\
		\end{pmatrix}.\]
\qedhere
\end{proof}
	
\begin{remark}
	In the proof, we have explicitly described the homomorphism $\xi$ in \eqref{xi} of Section \ref{prelim} in the particular case where $T$ is the norm torus.
\end{remark}

Recall from \eqref{iso} that  $T[p^n] \simeq \left(\Z/p^n\Z\right)^{p-1}$, hence we have a natural identification of $\Hh^1_{\loc}\left(G, T[p^n]\right)$ with $\Hh^1_{\loc}\left(G,\left(\Z/p^n\Z\right)^{p-1}\right)$ by inducing the action of $G$ on $\left(\Z/p^n\Z\right)^{p-1}$ via the same isomorphism. In the following proposition we show that these groups are non-trivial.
	
\begin{proposition}\label{h1loc}
\ja{Let $p$ be an odd prime and let $n \geq 2$}. Consider the action of $G$ on $\left(\Z/p^n\Z\right)^{p-1}$ induced by the isomorphism $T[p^n] \simeq \left(\Z/p^n\Z\right)^{p-1}$ of \eqref{iso}.
	There exists a (unique) extension of
	\[
		\gamma_1 \longmapsto v_1 = \begin{pmatrix} p^{n-2}(p-1)\\0\\ \vdots\\0\\p^{n-2}\\\end{pmatrix}, \quad 
		\gamma_2 \longmapsto v_2 = \begin{pmatrix} p^{n-1} \\ \vdots\\p^{n-1}\\0\\\end{pmatrix}\\
	\]
	to a cocycle in $\Hh^1\left(G,\left(\Z/p^n\Z\right)^{p-1}\right)$ and it is a non-trivial element of $\Hh_{\loc}^1\left(G, \left(\Z/p^n\Z\right)^{p-1}\right)$.
\end{proposition}
\begin{proof}
	To check that the assigned vectors define a cocycle, we have to prove that in $\left(\Z/p^n\Z\right)^{p-1}$
	\begin{align}
		&(1+\gamma_1+\dots+\gamma_1^{p-1})v_1 = \underline{0} \label{gamma1}\\ 
		&(1+\gamma_2+\dots+\gamma_2^{p^{n-1}-1})v_2 = \underline{0} \label{gamma2} \\ 
		&(1-\gamma_2)v_1 + (\gamma_1-1)v_2 = \underline{0} 	 \label{gamma12}
	\end{align}
	where $\underline{0}$ is the vector with all coordinates equal to zero; indeed these conditions must be true since $\gamma_1^p=\gamma_2^{p^{n-1}}=\gamma_1\gamma_2\gamma_1^{-1}\gamma_2^{-1} = 1$ in $G$.
	A lifting of $\gamma_1$ to $\GL_{p-1}(\Q)$ is the matrix $\gamma_1$ itself. It solves the polynomial $x^p-1$ and it is easy to see that $1$ is not an eigenvalue. Thus the minimal polynomial of $\gamma_1$ is $x^{p-1}+x^{p-2} + \ldots + x+1$ and so \eqref{gamma1} holds. Using $\gamma_2 = (p+1)\Id$, we have that
    \[(1+\gamma_2+\ldots+\gamma_2^{p-1})v_2 = (1 +(p+1) + \ldots + (p+1)^{p-1})v_2 \equiv \underline{0} \mod p^n.\]
    Since $n \geq 2$, we can collect the factor $1 +\gamma_2 + \ldots + \gamma_2^{p-1}$ on the left-hand side of \eqref{gamma2}, hence    \eqref{gamma2} holds.
    For \eqref{gamma12}, some simple calculations lead to:
	\[
		(1-\gamma_2)v_1 + (\gamma_1-1)v_2 \equiv \begin{pmatrix} p^{n-1} \\ 0 \\ \vdots \\ 0 \\ -p^{n-1} \\ \end{pmatrix} + \begin{pmatrix} -p^{n-1} \\ 0 \\ \vdots \\ 0 \\ p^{n-1}\\ \end{pmatrix} \equiv \underline{0} \mod p^n.
	\]
	Thus we can extend $v_1$ and $v_2$ to a cocycle $Z = \{Z_{\gamma}\}_{\gamma\in G}$, with $Z_{\gamma_1}=v_1$ and $Z_{\gamma_2} = v_2$. We now prove that $Z$ is not trivial. Suppose that it is a coboundary, then there exists $w \in \left(\Z/p^n\Z\right)^{p-1}$ such that $v_1 = (\gamma_1-1)w$ and $v_2 = (\gamma_2-1)w$. We denote with $w^{(i)}$ the $i$--th coordinate of $w$. From $v_2 = (\gamma_2-1)w$ we have $v_2 = p w$, so that $pw^{(p-1)} = 0$. On the other hand, from $v_1 = (\gamma_1-1)w$ we have	
	\begin{equation}\label{condition1}
		v_1 =
		\begin{pmatrix}
			p^{n-2}(p-1)\\
			0\\
			\vdots\\
			0\\
			p^{n-2}\\
		\end{pmatrix} = 
		\begin{pmatrix}
			-1 &  &  &  &  -1\\
			1 & -1 &  & & -1\\
			& \ddots & \ddots & &\vdots\\
			&  & 1 &-1 & -1\\
			&  &  & 1 & -2 \\
		\end{pmatrix}
		\begin{pmatrix}
			w^{(1)}\\
			w^{(2)}\\
			\vdots\\
			\vdots\\
			w^{(p-1)}\\
		\end{pmatrix}.
	\end{equation}
	Thus we obtain the system
	\begin{equation}\label{system}
		\begin{cases*}
			-w^{(1)} - w^{(p-1)} = p^{n-2}(p-1)\\
			w^{(1)} - w^{(2)} - w^{(p-1)} = 0\\
			\vdots\\
			w^{(p-3)}-w^{(p-2)}-w^{(p-1)} = 0\\
			w^{(p-2)}-2w^{(p-1)} = p^{n-2}.\\
		\end{cases*} 
	\end{equation}
	By adding the equations in \eqref{system} we get $p^{n-1} = -pw^{(p-1)}$, in contradiction with $pw^{(p-1)} = 0$. Therefore $Z$ is not a trivial cocycle in $\Hh^1\left(G,\left(\Z/p^n\Z\right)^{p-1} \right)$.\\
	We are left to show that $Z$ satisfies the local conditions. It is easy to see that the elements $\gamma_2^h$ and $\gamma_1\gamma_2^h$, for $h = 0,1,\dots,p^{n-1}-1$ are generators of all the cyclic subgroups of $G$. So it is enough to show that there exist $W_{\gamma_2^h}$ and $W_{\gamma_1\gamma_2^h}$ in $(\Z/p^n\Z)^{p-1}$ such that $Z_{\gamma_2^h} = (\gamma_2^h-1) W_{\gamma_2^h}$ and $Z_{\gamma_1\gamma_2^h} = (\gamma_1\gamma_2^h - 1)W_{\gamma_1\gamma_2^h}$ for all $h \in \{0,1,\dots,p^{n-1}-1\}$.
	First, since $\gamma_2-1 = p\Id$, we have that the image of $\gamma_2-1$ is the submodule $M$ of $(\Z/p^n\Z)^{p-1}$ of vectors with each coordinate divisible by $p$. The vector $v_2$ satisfies this condition, so there exists $W_{\gamma_2} \in \left(\Z/p^n\Z\right)^{p-1}$ such that $v_2 = Z_{\gamma_2} = (\gamma_2-1)W_{\gamma_2}$.
    Furthermore, since $Z_{\gamma_2^h} = (1 + \gamma_2 + \ldots + \gamma_2^{h-1})v_2$, we have that $Z_{\gamma_2^h}$ also lies in $M$ and $Z_{\gamma_2^h} = (\gamma_2^h-1)W_{\gamma_2}$.
	Now we fix $h \in \{0,1,\dots, p^{n-1}-1\}$ and define 
	\[V = \left\{ \begin{pmatrix} v^{(1)}\\ \vdots\\ v^{(p-1)}\\ \end{pmatrix} \in \left(\Z/p^n\Z\right)^{p-1}\, \Bigg\vert \, \sum_{j = 1}^{p-1} v^{(j)} \equiv 0 \mod p \right\} \subseteq \left(\Z/p^n\Z\right)^{p-1}.\]
	We claim that the image of $\gamma_1\gamma_2^h - 1$ is equal to $V$.
	Since $\gamma_2^h = (p+1)^h \Id \equiv (1+pl) \Id \mod p^n$ for some $l \in \Z/p^n\Z$, we can rewrite
	\[\gamma_1\gamma_2^h -1 = (1+pl)\gamma_1 -1=
	   \begin{pmatrix} 
			-1 &  & &   & -1-pl\\
			1+pl & -1 & & &-1-pl\\
			& \ddots & \ddots & & \vdots\\
			& & \ddots & -1& -1-pl\\
			& &    & 1+pl &-2-pl\\
		\end{pmatrix}
	\]
	\noindent and thus $\imm(\gamma_1\gamma_2^h-1) \subseteq V$.\\
	With easy calculations one can check that the determinant of $\gamma_1\gamma_2^h-1$ is equal to $p$ modulo $p^2$. 
    Take a lifting of $\gamma_1\gamma_2^h-1$ to an integer matrix; this integer matrix has still determinant equal to $p$ modulo $p^2$. Since $\Z$ is a PID, we can consider its Smith normal form $\diag(\alpha_1,\alpha_2,\ldots,\alpha_{p-1})$, and we get that $p \mid \alpha_{p-1}$ while $p^2 \nmid \alpha_{p-1}$ and $p \nmid \alpha_j$ for every $j =1, \ldots, p-2$.
    Its projection $\diag(\wt \alpha_1,\wt \alpha_2, \ldots, \wt \alpha_{p-1})$ modulo $p^n$ is such that $\wt \alpha_j$ is invertible, for $j = 1,\ldots, p-2$, and $\wt \alpha_{p-1}\neq 0$ is equal to $0$ modulo $p$. Therefore, up to basis changes in $(\Z/p^n\Z)^{p-1}$, the map $\gamma_1\gamma_2^h-1$ is
	\[
	   \begin{pmatrix}
			1&  &  & \\
			& \ddots & & \\
			& & 1& \\
			& & & p\\
		\end{pmatrix}.
	\]
    It follows that $\imm(\gamma_1\gamma_2^h-1)$ has index equal to $p$ in $\left(\Z/p^n\Z\right)^{p-1}$. The submodule $V$ has also index equal to $p$ in $\left(\Z/p^n\Z\right)^{p-1}$ and so, from the inclusions $\imm(\gamma_1\gamma_2^h-1) \subseteq V \subseteq \left(\Z/p^n\Z\right)^{p-1}$, we get the equality $V = \imm(\gamma_1\gamma_2^h-1)$.
	To conclude that $Z$ satisfies the local conditions it only remains to verify that $Z_{\gamma_1\gamma_2^h}$ lies in $V$. We have $Z_{\gamma_1\gamma_2^h} = v_1 + \gamma_1Z_{\gamma_2^h}$, with $v_1 \in V$ and, as mentioned above, $Z_{\gamma_2^h} \in M$ (and also $\gamma_1Z_{\gamma_2^h} \in M$). 
    Since $M$ is contained in $V$, we get that $Z_{\gamma_1\gamma_2^h} \in V$.
	\qedhere
\end{proof}
We conclude this section with the proof of Theorem \ref{thm:1}.
\begin{proof}[Proof of Theorem \ref{thm:1}.]
    Let $p$ be an odd prime. 
	\par (a) Let $T$ be an algebraic torus of dimension $r<p-1$ defined over a number field $k$ and let $n \geq 1$ be an integer. With the notation of Section \ref{prelim}, since the inclusion $K\subseteq L\left(\z\right)$ holds, it is clear that $G$ is isomorphic to the quotient $\Gal(L\left(\z\right)/k)/\Gal(L\left(\z\right)/K)$. Since $p$ is an odd prime, any $p$--Sylow subgroup of the group $\Delta \simeq \Gal(L/k)$ is contained in $\SL_r(\Q)$; hence by the condition $r<p-1$ and Lemma \ref{illengo}, we have that $\Delta$ has no nontrivial $p$--Sylow subgroups. Since
	\[ \Gal(L(\z)/k) \longhookrightarrow \Gal(L/k) \times \Gal(k(\z)/k),
	\]
	we get that any $p$--Sylow subgroup of $\Gal(L(\z)/k)$ is isomorphic to a subgroup of $\Gal(k(\z)/k)$. The latter is isomorphic to a subgroup of $\left(\Z/p^n\Z\right)^{\times}$, hence it is cyclic. Thus $\Gal(L(\z)/k)$ contains only one cyclic $p$--Sylow subgroup. Now let $G_p$ be a $p$--Sylow subgroup of $G$. By Theorem \ref{sylow quot}, it is the image through the projection to the quotient of the $p$--Sylow subgroup of $\Gal(L(\z)/k)$, hence it is cyclic too and $\Hh^1_{\loc}(G_p, T[p^n])$ is trivial.

	\par (b) Let $n\geq 2$ be an integer. As mentioned in the Introduction, it is enough to prove the statement for $r = p-1$. By Lemma \ref{costr toro} and Proposition \ref{h1loc} there exists an algebraic torus defined over $k = \Q(\z_p)$ of dimension $p-1$ such that $\Gal(k(T[p^n]/k)$ is isomorphic to $\Z/p\Z\times \Z/p^{n-1}\Z$ and $\Hh^1_{\loc}(G,T[p^n]) \neq 0$. By Theorem \ref{neccond}, there exists a number field $L$ such that $L\cap k\left(T[p^n]\right) = k$ and the local-global divisibility by $p^n$ does not hold for $T(L)$.
	\qedhere
\end{proof}

\begin{remark} \label{rational_counterexamples}
	The example built in Lemma \ref{costr toro} is defined over the number field $k = \Q(\z_p)$, but we can use it to construct an example over $\Q$ of dimension $(p-1)^2$. Indeed, let $\wt T = R_{k/\Q}(T)$ be the Weil restriction of $T$. It has dimension $\dim(\wt T) = [k:\Q] \dim(T) = (p-1)^2$ and for every number field $F$ containing $k$ we have $\wt T(F) = R_{k/\Q}(T)(F) \simeq T(F)^{p-1}$. In particular, if $P\in T(F)$ is a point such that the local-global divisibility fails, we have the failure  of the the local-global divisibility also for the corresponding point on $\wt T(F)$ (given by $p-1$ copies of $P$).
\end{remark}

\lp{

\begin{remark} \label{extended_counterexamples}
Let $T$ be a torus defined over a number field $k$, with non-trivial $\Hh^1_\loc\left(\Gal(k(T[p^n])/k),T[p^n]\right)$.
If $L$ is a finite extension of $k$ linearly disjoint from $k(T[p^n])$ over $k$,
then
$$\Hh^1_{\loc}\left(\Gal(L(T[p^n])/L),T[p^n]\right) \simeq \Hh^1_{\loc}\left(\Gal(k(T[p^n])/k),T[p^n]\right)$$
is non-trivial too.
Thus by Theorem \ref{neccond} we have a counterexample over a finite extension of $L$.
In this way, we have counterexamples over infinitely many number fields; in particular this applies to $k = \Q$.
Moreover, this argument works not only for tori, but for every commutative algebraic group.
\end{remark}

}

\section{Proof of Theorem \ref{thm:2}} \label{sec:dim2}	

We are going to prove Theorem \ref{thm:2}. In the following proof, we will show that $\Hh^1_{\loc}(G, T[p^n]) = 0$; in particular this proves that the local-global divisibility by $p^n$ holds.
	
\begin{proof}[Proof of Theorem \ref{thm:2}.]
We proceed by induction on $n \geq 1$ proving that $\Hh^1_{\loc}(G, T[p^n]) = 0$. The base of the induction $n=1$ is proven in \cite{Ill}. Thus suppose $n\geq 2$ and $\Hh^1_{\loc}(G,T[p^m]) = 0$, for every $m\leq n-1$.
With the notation of Section \ref{prelim}, let $F = K \cap k(\z)$. We have that $[F:k] = [G:G']$, thus it is coprime with $p$ by the assumptions of the theorem.
Therefore the $p$--Sylow subgroups of $G$ are contained in $G'$, hence in $\wt\Delta$. 
For any $j = 1,\dots,n-1$, using the same construction as in Section \ref{prelim}, we can define a group $G^{(j)} \subseteq \GL_r(\Z/p^j\Z)$ that is isomorphic to $\Gal(k\left(T[p^j]\right)/k)$.
Let $\pi_j: \GL_r(\Z/p^n\Z) \longrightarrow \GL_r(\Z/p^j\Z)$ be the reduction modulo $p^j$ and let $H^{(j)}$ be the intersection of $G$ with $\ker \pi_j$. It is easy to prove that $G^{(j)}= \pi_j(G) \simeq \Gal(k\left(T[p^j]\right)/k)$ and $H^{(j)} \simeq \Gal(K/k\left(T[p^j]\right))$. We claim that $H^{(1)}$ is trivial. Every $h\in H^{(1)}$ can be written as $h = 1 + pA$, for some $A \in \Mat_r(\Z/p^n\Z)$, so we have $h^{p^{n-1}} \equiv \Id \mod \it{p^n}$. Hence the subgroup $H^{(1)}$ is a $p$--group and thus it is contained in a $p$--Sylow subgroup $G_p$ of $G$.
Since every $p$--Sylow subgroup of $G$ is contained in $G'$, we have that $H^{(1)}$ is contained in a $p$--Sylow subgroup $G'_p$ of $G'$ too. Thus $H^{(1)}\subseteq \wt\Delta_p$, where $\wt\Delta_p$ is a $p$--Sylow subgroup of $\wt\Delta$. Let $H \subset \GL_r(\Z)$ be a subgroup of $\Delta$ such that its image via the reduction modulo $p^n$ is $H^{(1)}$. Consider $\pi:\GL_r(\Z) \longrightarrow \GL_r(\Z/p\Z)$ the reduction modulo $p$. We have that $H$ is contained in $\ker(\pi) \cap \Delta$. So, by using Lemma \ref{lemma:inj} we find that $H$ and, consequently, $H^{(1)}$ are trivial. Our claim is proved and we get that also every $H^{(j)}$ is trivial since they are all contained in $H^{(1)}$. Therefore $G \simeq G^{(n-1)} \simeq \dots \simeq G^{(1)}$ via the relative projections. 

Consider the following short exact sequence
\[ 1 \longrightarrow T[p]\stackrel{\iota}{\longrightarrow} T[p^n] \stackrel{\varepsilon}{\longrightarrow} T[p^{n-1}] \longrightarrow 1, \]
where $\iota$ is the inclusion and $\varepsilon$ is the $p$--power map (here we are using the multiplicative notation for $T[p]$, $T[p^{n-1}]$ and $T[p^n]$). 
The group $G$ acts on $T[p^n]$ and, via the projections, on $T[p]$ and $T[p^{n-1}]$ and by these actions the above short exact sequence is a sequence of $G$--modules. Thus we have the following long exact sequence:
\[ 1 \rightarrow T[p]^{G} \rightarrow T[p^n]^{G} \rightarrow T[p^{n-1}]^{G} \rightarrow \Hh^1(G,T[p]) \rightarrow \Hh^1(G,T[p^n]) \rightarrow \Hh^1(G,T[p^{n-1}]) \rightarrow \dots \]
\noindent Let $C$ be a cyclic subgroup of $G$ and for $i= 1, n-1, n$ let $\res_i$ be the restriction $\Hh^1(G,T[p^i])\longrightarrow \Hh^1(C,T[p^i])$. We have the following diagram with exact rows
\begin{figure}[H]
	\centering
	\begin{tikzpicture}\color{black}
		\node(a) at (-3,2) {$\ker(\res_1)$};
		\node(b) at (0,2) {$\ker(\res_n)$};
		\node(c) at (3,2) {$\ker(\res_{n-1})$};
		\node (d) at (-3,0) {$\Hh^1(G,T[p])$};
		\node(e) at (0,0) {$\Hh^1(G,T[p^n])$};
		\node(f) at (3,0) {$\Hh^1(G,T[p^{n-1}])$};
		\node (g) at (-3,-2) {$\Hh^1(C,T[p])$};
		\node (h) at (0,-2) {$\Hh^1(C,T[p^n])$};
		\node(i) at (3,-2) {$\Hh^1(C,T[p^{n-1}])$};
			
		\node at (-2.3,-1) {$\res_1$};
		\node at (0.5,-1) {$\res_n$};
		\node at (3.7, -1) {$\res_{n-1}$};
			
		\draw[-stealth] (a) --(b);
		\draw[-stealth] (b) --(c);
		\draw[-stealth] (a) --(d);
		\draw[-stealth] (b) --(e);
		\draw[-stealth] (c) --(f);
		\draw[-stealth] (d) --(e);
		\draw[-stealth] (e) --(f);
		\draw[-stealth] (d) --(g);
		\draw[-stealth] (e) --(h);
		\draw[-stealth] (f) --(i);
		\draw[-stealth] (g) --(h);
		\draw[-stealth] (h) --(i);
	\end{tikzpicture}
\end{figure}
\noindent where the central row is given by the long exact sequence above, the lower row is obtained from the same exact sequence by restriction to the subgroup $C$, while the upper one is induced by the commutativity of the diagram given by the last two rows.
Since the group $G$ is isomorphic to $G^{(1)}$ and to $G^{(n-1)}$, we have that the cyclic subgroups of $G^{(1)}$ and $G^{(n-1)}$ are the images of the projections of the cyclic subgroups of $G$. Therefore, by taking the  intersection over all the cyclic subgroups of $G$, from the first row of the diagram we have the exact sequence
\[ \Hh^1_{\loc}(G,T[p]) \longrightarrow \Hh^1_{\loc}(G,T[p^n]) \longrightarrow \Hh^1_{\loc}(G,T[p^{n-1}]). \]
Since by inductive hypothesis $ \Hh^1_\loc(G,T[p]) = \Hh^1_{\loc}(G,T[p^{n-1}])=0$, we also have $\Hh^1_{\loc}(G,T[p^n]) =0$.

\qedhere
\end{proof}
	
\begin{remark}
	We remark that in Lemma \ref{costr toro} we have $k = \Q(\z_p)$ and $F := k\left(T[p^n]\right) \cap k(\z_{p^n}) = k(\z_{p^n})$, so $[F:k] = [k(\z_{p^n}):k] = [\Q(\z_{p^n}) : \Q(\z_p)]=p^{n-1}$. Moreover, we notice that the automorphism $\eta$ in the proof of Lemma \ref{costr toro} is a generator of the subgroup $H^{(1)} \simeq \Gal\left(k\left(T[p^n]\right)/k\left(T[p]\right)\right)$, defined in the proof of Theorem \ref{thm:2} above. In particular, for the torus defined in Lemma \ref{costr toro} we have that $H^{(1)}$ is not trivial.
\end{remark}

\begin{remark}
	Notice that it  is always possible to construct an algebraic torus $T$, not split over $k$, that satisfies the conditions of Theorem \ref{thm:2}. An example is obtained by taking a number field $k$ that contains a $p^n$--th root of unity, any finite extension $L/k$ of degree $d$, with $p-1 \leq d <3(p-1)$, and considering $T = R_{L/k}\G_{m,L}$, the Weil restriction of the split torus $\G_{m,L}$ over $L$. The torus $T$ is defined over $k$, it is split over $L$ and it has dimension $d$. Thus $F = K\cap k(\z)$ is equal to $k$ and the hypotheses of the theorem are satisfied.
\end{remark}

\bigskip
\bigskip

\begin{minipage}[t]{10cm}
	\begin{flushleft}
		\small{
			\textsc{Jessica Alessandr\`i}
			\\* Universit\`a degli Studi dell’Aquila,
			\\* Via Vetoio, Coppito 1
			\\* Coppito (AQ), 67100, Italy
			\\*e-mail: jessica.alessandri@graduate.univaq.it
				
		}
	\end{flushleft}
\end{minipage}
	
\bigskip

\begin{minipage}[t]{10cm}
    \begin{flushleft}
	   	\small{
	    	\textsc{Rocco Chiriv\`i}
		    \\* Universit\`a del Salento,
		    \\* Via per Arnesano
		    \\* Monteroni di Lecce (LE), 73047, Italy
		    \\*e-mail: rocco.chirivi@unisalento.it
			
	    }
	   \end{flushleft}
\end{minipage}

\bigskip
	
\begin{minipage}[t]{10cm}
    \begin{flushleft}
	    \small{
		    \textsc{Laura Paladino}
   		    \\*Universit\`a della Calabria,
   		   	\\* Ponte Bucci, Cubo 30B 
	       	\\* Rende (CS), 87036, Italy
	    	\\*e-mail: laura.paladino@unical.it
			
    	}
    \end{flushleft}
\end{minipage}

\end{document}